\documentclass[12pt]{amsart}
\usepackage[cp1251]{inputenc}
\usepackage[T2A]{fontenc}
\usepackage[english]{babel}
\usepackage{amsmath,amsfonts,amssymb}
\usepackage{geometry}
\usepackage{amsmath, amsthm, amscd, amsfonts, amssymb, graphicx, color}
\usepackage[bookmarksnumbered, colorlinks, plainpages]{hyperref}

\theoremstyle{remark}

\theoremstyle{definition}

\begin{document}

\title[ELLIPTIC BOUNDARY-VALUE PROBLEMS IN THE SENSE OF LAWRUK]{ELLIPTIC BOUNDARY-VALUE PROBLEMS\\IN THE SENSE OF LAWRUK ON\\ SOBOLEV AND H\"ORMANDER SPACES}

\author[I. Chepurukhina]{Iryna S. Chepurukhina}

\address{Institute of Mathematics, National Academy of Sciences of Ukraine,
3 Tereshchenkivs'ka, Kyiv, 01601, Ukraine}

\email{Chepuruhina@mail.ru}


\author[A. Murach]{Aleksandr A. Murach}

\address{Institute of Mathematics, National Academy of Sciences of Ukraine,
3 Tereshchenkivs'ka, Kyiv, 01601, Ukraine}

\email{murach@imath.kiev.ua}



\begin{abstract}
We investigate elliptic boundary--value problems with additional unknown functions in boundary conditions. These problems were introduced by Lawruk. We prove that the operator corresponding to such a problem is bounded and Fredholm on appropriate couples of the inner product isotropic H\"ormander spaces $H^{s,\varphi}$, which form the refined Sobolev scale. The order of differentiation for these spaces is given by the real number $s$ and positive function $\varphi$ that varies slowly at infinity in the sense of Karamata. We consider this problem for an arbitrary elliptic equation $Au=f$ on a bounded Euclidean domain $\Omega$ under the condition that $u\in H^{s,\varphi}(\Omega)$,  $s<\mathrm{ord}\,A$, and $f\in L_{2}(\Omega)$. We prove theorems on the a priori estimate and regularity of the generalized solutions to this problem.
\end{abstract}

\maketitle

\section{Introduction}
In this paper, we investigate
elliptic boundary-value problems with additional unknown functions in boundary conditions. These functions are given on the smooth boundary of a Euclidean domain $\Omega$, where the problem is posed. Such problems were introduced by B.~Lawruk \cite{Lawruk63a, Lawruk63b, Lawruk65}. They appear naturally if we pass from a general (nonregular) elliptic  boundary-value problem to its formally adjoint problem. Moreover, the class of these problems are closed with respect to this passage. Their various examples occur in hydrodynamics and the theory of elasticity \cite{AslanyanVassilievLidskii81, Ciarlet90, NazarovPileckas93}.

This class has been investigated in some function spaces by V.~A. Kozlov, V.~G. Maz'ya, and J.~Rossmann \cite{KozlovMazyaRossmann97} (Chapt.~3), mainly, for  scalar elliptic equations and by I.~Ya.~Roitberg \cite{RoitbergInna97, RoitbergInna98} for mixed-order elliptic systems. The results of I.~Ya.~Roitberg is expounded in Ya.~A.~Roitberg's monograph \cite[Chapt.~2]{Roitberg99}. The main result obtained are the theorem on the Fredholm property of these problems on two-sided scales of normed spaces, solutions of the elliptic equations being considered in the spaces introduced by Ya.~A.~Roitberg \cite{Roitberg64, Roitberg65} (see also his monograph \cite{Roitberg96} (Sect.~2)). These spaces coincide with the Sobolev spaces provided that their differentiation order is large enough. But, for other orders, the Roitberg spaces contain elements that are not distributions on $\Omega$.

The purpose of this paper is to prove a version of this theorem for Sobolev and more general H\"ormander spaces \cite{Hermander63} (Sect.~2.2) of distributions on $\Omega$. We use the H\"ormander spaces
$$
H^{s,\varphi}(\mathbb{R}^{n}):=
\bigl\{w\in\mathcal{S}'(\mathbb{R}^{n}):\,
\langle\xi\rangle^{s}\varphi(\langle\xi\rangle)\widehat{w}(\xi)\in
L_{2}(\mathbb{R}^{n},d\xi)\bigr\}
$$
and their analogs for the domain $\Omega$ and its boundary. These spaces are  parametrized with the number $s\in\mathbb{R}$ and the function $\varphi:[1,\infty)\rightarrow(0,\infty)$ that varies slowly at infinity in the sense of J.~Karamata. Here, $\widehat{w}$ is the Fourier transform of the tempered distribution $w$, whereas $\langle\xi\rangle:=(1+|\xi|^{2})^{1/2}$. If $\varphi\equiv1$, then $H^{s,\varphi}(\mathbb{R}^{n})$ becomes the Sobolev space $H^{s}(\mathbb{R}^{n})$.

These H\"ormander spaces form the refined Sobolev scale, which was selected and investigated by V.~A.~Mikhailets and A.~A.~Murach \cite{MikhailetsMurach05UMJ5, MikhailetsMurach06UMJ3,  MikhailetsMurach08MFAT1}. This scale has an important interpolation property; namely, every space $H^{s,\varphi}(\mathbb{R}^{n})$ is a result of the interpolation with an appropriate function parameter of the inner product Sobolev spaces $H^{s-\varepsilon}(\mathbb{R}^{n})$ and $H^{s+\delta}(\mathbb{R}^{n})$ with $\varepsilon,\delta>0$.

V.~A.~Mikhailets and A.~A.~Murach [15, 16, 18~-- 23] have elaborated the theory of solvability of general elliptic boundary-value problems on the refined Sobolev scale and its modification by Ya.~A.~Roitberg. The above-mentioned interpolation is a key method in their theory. In this connection, we also note papers \cite{Slenzak74, AnopMurach14MFAT2, AnopMurach14UMJ7}. However, this theory does not involve the important class of elliptic boundary-value problems in the sense of B.~Lawruk.

In this paper, we consider these problems for an arbitrary scalar   elliptic equation $Au=f$ with $u\in H^{s,\varphi}(\Omega)$ for $s<\mathrm{ord}\,A$ and $f\in L_{2}(\Omega)$. This approach originates from the papers by J.-L.~Lions and E.~Magenes \cite{LionsMagenes62, LionsMagenes63} (see also \cite{Murach09MFAT2}). The case of $s\geq\mathrm{ord}\,A$ was investigated in~\cite{Chepurukhina14Coll2}.

This paper consists of seven sections. Section~1 is Introduction. In Section~2, we formulate an elliptic boundary-value problem in the sense of B.~Lawruk and discuss the corresponding Green formula and formally adjoint problem. In Section~3, we give the definitions of H\"ormander function spaces that form the refined Sobolev scales. Section~4 contains the main results of the paper. They are the theorem on the Fredholm property of the problem under investigation, the a priory estimate for its generalized solutions, and the theorem about their regularity in the refined Sobolev scales. In Sections 5 and 6, we discuss auxiliary results that we need to prove these theorems. The proofs of the main results are given in Section~7.

\section{Statement of the problem} 

Let $\Omega$ be a bounded domain in $\mathbb{R}^{n}$, where $n\geq2$, with the infinitely smooth boundary $\Gamma:=\partial\Omega$, which is a closed manifold of dimension $n-1$. Let $\nu(x)$ denote the unit vector of the inward normal to $\Gamma$ at a point $x\in\Gamma$.

Choose arbitrarily integers $q\geq1$, $\varkappa\geq1$, $m_{1},\ldots,m_{q+\varkappa}\in[0,2q-1]$ and $r_{1},\ldots,r_{\varkappa}$.
In the paper, we consider the linear boundary-value problem on the domain $\Omega$ with $\varkappa$ additional unknown functions on the boundary $\Gamma$:
\begin{gather}\label{3f1}
Au=f\quad\mbox{on}\quad\Omega,\\
B_{j}\,u+\sum_{k=1}^{\varkappa}C_{j,k}\,v_{k}=g_{j}\quad\mbox{on}\quad\Gamma,
\quad j=1,...,q+\varkappa.\label{3f2}
\end{gather}
Here,
$$
A:=A(x,D):=\sum_{|\mu|\leq 2q}a_{\mu}(x)D^{\mu}
$$
is a linear differential operator on $\overline{\Omega}:=\Omega\cup\Gamma$ of the even order $2q$. Besides, each
$$
B_{j}(x,D):=\sum_{|\mu|\leq m_{j}}b_{j,\mu}(x)D^{\mu}
$$
is a linear boundary differential operator on $\Gamma$ of order $m_{j}$, and every $C_{j,k}:=C_{j,k}(x,D_{\tau})$ is a linear tangent differential operator on $\Gamma$ with $\mathrm{ord}\,C_{j,k}\leq m_{j}+r_{k}$. All the coefficients of these differential operators are supposed to be infinitely smooth complex-valued functions given on $\overline{\Omega}$ and $\Gamma$ respectively.

Here and below, we use the following standard notation:
$\mu:=(\mu_{1},\ldots,\mu_{n})$ is a multi-index with $|\mu|:=\mu_{1}+\ldots+\mu_{n}$, and
$D^{\mu}:=D_{1}^{\mu_{1}}\ldots D_{n}^{\mu_{n}}$ with  $D_{\ell}:=i\partial/\partial x_{\ell}$, where $i$ is imaginary unit and  $x=(x_1,\ldots,x_n)\in\mathbb{R}^{n}$. Besides, $D_{\nu}:=i\partial/\partial\nu(x)$ and  $\xi^{\mu}:=\xi_{1}^{\mu_{1}}\ldots\xi_{n}^{\mu_{n}}$ for  $\xi=(\xi_{1},\ldots\xi_{n})\in\mathbb{C}^{n}$.

The function $u$ given on $\Omega$ and all the functions  $v_{1},\ldots,v_{\varkappa}$ given on $\Gamma$ are unknown in the
problem \eqref{3f1}, \eqref{3f2}. In the paper, all functions and distributions are supposed to be complex-valued.

We assume that the boundary-value problem \eqref{3f1}, \eqref{3f2} is elliptic on $\Omega$ in the sense of B.~Lawruk~\cite{Lawruk63a}. Recall the corresponding definition (see also \cite{KozlovMazyaRossmann97} (Subsect.~3.1.3)).

Let $A^{(0)}(x,\xi)$, $B_{j}^{(0)}(x,\xi)$, and $C_{j,k}^{(0)}(x,\tau)$
denote the principal symbols of the differential operators $A(x,D)$, $B_{j}(x,D)$, and $C_{j,k}(x,D_{\tau})$ respectively. We recall that
$$
A^{(0)}(x,\xi):=\sum_{|\mu|=2q}a_{\mu}(x)\xi^{\mu},\quad\mbox{with}\;\; x\in\overline{\Omega}\;\;\mbox{and}\;\;\xi\in\mathbb{C}^{n},
$$
is a homogeneous polynomial of order $2q$ in $\xi$ and that
$$
B_{j}^{(0)}(x,\xi):=\sum_{|\mu|=m_{j}}b_{j,\mu}(x)\xi^{\mu},\quad\mbox{with}\;\; x\in\Gamma\;\;\mbox{and}\;\;\xi\in\mathbb{C}^{n},
$$
is a homogeneous polynomial of order $m_{j}$ in $\xi$. Besides, for every point $x\in\Gamma$, the expression $C^{(0)}_{j,k}(x,\tau)$ is a homogeneous polynomial of order $m_{j}+r_{k}$ in $\tau$, with $\tau$ being a tangent vector to the boundary $\Gamma$ at the point $x$. (Of course, if $\mathrm{ord}\,C_{j,k}<m_{j}+r_{k}$, then $C^{(0)}_{j,k}(x,\tau):=\nobreak0$.)

The boundary-value problem \eqref{3f1}, \eqref{3f2} is said to be elliptic on $\Omega$ if the following three conditions are fulfilled:

(i) The differential operator $A(x,D)$ is elliptic at each point $x\in\overline{\Omega}$; i.e., $A^{(0)}(x,\xi)\neq\nobreak0$ for an arbitrary vector $\xi\in\mathbb{R}^{n}\setminus\{0\}$.

(ii) The differential operator $A(x,D)$ is properly elliptic at each point $x\in\Gamma$; i.e., for an arbitrary vector $\tau\neq0$ tangent to $\Gamma$ at $x$, the polynomial $A^{(0)}(x,\tau+\zeta\nu(x))$ in $\zeta\in\mathbb{C}$ has $q$ roots with positive imaginary part and $q$ roots with negative  imaginary part (these roots are calculated with regard for their multiplicity).

(iii) The system of boundary-value conditions \eqref{3f2} covers equation \eqref{3f1} at each point $x\in\Gamma$. This means that, for an arbitrary vector $\tau\neq0$ from condition (ii), the boundary-value problem
\begin{gather*}
A^{(0)}(x,\tau+D_{t}\,\nu(x))\theta(t)=0\quad\mbox{for}\quad t>0,\\
B_{j}^{(0)}(x,\tau+D_{t}\,\nu(x))\theta(t)\big|_{t=0}+
\sum_{k=1}^{\varkappa}C_{j,k}^{(0)}(x,\tau)\lambda_{k}=0,\quad j=1,...,q+\varkappa,
\end{gather*}
has only zero solution. This problem is considered relative to the unknown function $\theta\in C^{\infty}([0,\infty))$ with $\theta(t)\rightarrow0$ as $t\rightarrow\infty$ and the unknown complex-valued numbers $\lambda_{1},\ldots,\lambda_{\varkappa}$. Here, $A^{(0)}(x,\tau+D_{t}\,\nu(x))$ and $B_{j}^{(0)}(x,\tau+D_{t}\,\nu(x))$ are differential operators with respect to $D_{t}:=i\partial/\partial t$. We obtain them if we put $\zeta:=D_{t}$ in the polynomials $A^{(0)}(x,\tau+\zeta\nu(x))$ and $B_{j}^{(0)}(x,\tau+\zeta\nu(x))$ in $\zeta$.

It is worthwhile to note that condition (ii) follows from condition (i) in the $n\geq3$ case.

Examples of elliptic boundary-value problems in the sense of B.~Lawruk are discussed in \cite{KozlovMazyaRossmann97} (Subsect.~3.1.5). Among them, we mention the following simple example
\begin{equation*}
\Delta u=f\;\;\mbox{on}\;\;\Omega,\quad
u+v=g_{1}\;\;\mbox{and}\;\;D_{\nu}u+D_{\tau}v=g_{2}\;\;\mbox{on}\;\;\Gamma.
\end{equation*}
Here, $n=2$, $\varkappa=1$, $\Delta$ is the Laplace operator, and  $D_{\tau}:=i\partial/\partial\tau$, with $\partial/\partial\tau$ being the derivative along the curve $\Gamma$. It is easy to check that this boundary-value problem is elliptic on $\Omega$.

To the problem \eqref{3f1}, \eqref{3f2} corresponds the linear mapping
\begin{equation}\label{3f3}
\begin{gathered}
\Lambda:\,(u,v_{1},...,v_{\varkappa})\mapsto
\biggl(Au,\,B_{1}u+\sum_{k=1}^{\varkappa}C_{1,k}\,v_{k},...,
B_{q+\varkappa}\,u+\sum_{k=1}^{\varkappa}C_{q+\varkappa,k}\,v_{k}\biggr),\\
\mbox{with}\quad u\in C^{\infty}(\overline{\Omega})\quad\mbox{end}\quad
v_{1},\ldots,v_{\varkappa}\in C^{\infty}(\Gamma).
\end{gathered}
\end{equation}
We will investigate properties of the extension (by continuity) of this mapping in suitable couples of H\"ormander function spaces and, in particular, Sobolev spaces.

In order to describe the range of this extension, we need the following Green formula (see \cite{KozlovMazyaRossmann97} (Theorem~3.1.1)):
\begin{equation*}
\begin{gathered}
(Au,w)_{\Omega}+\sum_{j=1}^{q+\varkappa}\biggl(B_{j}\,u+
\sum_{k=1}^{\varkappa}C_{j,k}\,v_{k},h_{j}\biggr)_{\Gamma}=\\
=(u,A^{+}w)_{\Omega}+\sum_{j=1}^{2q}\biggl(D_{\nu}^{j-1}u,K_{j}\,w+
\sum_{k=1}^{q+\varkappa}Q_{k,j}^{+}\,h_{k}\biggr)_{\Gamma}+
\sum_{j=1}^{\varkappa}\biggl(v_{j},
\sum_{k=1}^{q+\varkappa}C_{k,j}^{+}\,h_{k}\biggr)_{\Gamma},
\end{gathered}
\end{equation*}
where $u,w\in C^{\infty}(\overline{\Omega})$, $v=(v_{1},\ldots,v_{\varkappa})\in(C^{\infty}(\Gamma))^{\varkappa}$, and  $h=(h_{1},\ldots,h_{q+\varkappa})\in(C^{\infty}(\Gamma))^{q+\varkappa}$ are arbitrary, whereas $(\cdot,\cdot)_{\Omega}$ and $(\cdot,\cdot)_{\Gamma}$ are the inner products in the Hilbert spaces $L_{2}(\Omega)$ and $L_{2}(\Gamma)$ of functions square integrable over $\Omega$ and $\Gamma$ respectively.

Here, $A^{+}$ stands for the differential operator that is formally adjoint to $A$ with respect to $(\cdot,\cdot)_{\Omega}$. Besides, $C_{k,j}^{+}$ and $Q_{k,j}^{+}$ denote the tangent differential operators that are formally adjoint to $C_{k,j}$ and $Q_{k,j}$ respectively relative to $(\cdot,\cdot)_{\Gamma}$, the tangent differential operators $Q_{k,j}$
appearing in the representation of the boundary differential operators $B_{j}$ in the form
$$
B_{j}(x,D)=
\sum_{k=1}^{2q}Q_{j,k}(x,D_{\tau})D_{\nu}^{k-1}.
$$
(If $k>m_{j}$, then $Q_{j,k}=0$). Finally, $K_{j}:=K_{j}(x,D)$ is a certain linear boundary differential operator on $\Gamma$, with $\mathrm{ord}\,K_{j}\leq2q-j$.

Taking into account the Green formula, we consider the following boundary-value problem on $\Omega$ with $q+\varkappa$ additional unknown function on $\Gamma$:
\begin{gather}\label{3f4}
A^{+}w=\omega\quad\mbox{in}\quad\Omega,\\
\label{3f5}
K_{j}\,w+\sum_{k=1}^{q+\varkappa}Q_{k,j}^{+}\,h_{k}=\chi_{j}\quad
\mbox{on}\quad\Gamma,\quad j=1,...,2q,\\
\label{3f6}
\sum_{k=1}^{q+\varkappa}C_{k,j}^{+}\,h_{k}=\chi_{2q+j}
\quad\mbox{on}\quad\Gamma,\quad j=1,...,\varkappa.
\end{gather}
This problem is formally adjoint to the problem \eqref{3f1}, \eqref{3f2} with respect to the Green formula. It is worthwhile to note that the problem \eqref{3f1}, \eqref{3f2} is elliptic if and only if the problem \eqref{3f4}, \eqref{3f5}, \eqref{3f6} is elliptic as well (see \cite{KozlovMazyaRossmann97} (Theorem~3.1.2)).

\section{H\"ormander spaces and refined Sobolev scales}
 
Here, we give the definitions of the H\"ormander inner product spaces that form the refined Sobolev scales over $\mathbb{R}^{n}$, $\Omega$, and $\Gamma$. We also discuss some of their properties \cite{MikhailetsMurach14}.

We should begin with the definition of the H\"ormander spaces $H^{s,\varphi}(\mathbb{R}^{n})$, with $s\in\mathbb{R}$ and $\varphi\in\mathcal{M}$. Here and below, $\mathcal{M}$ stands for the set of all Borel measurable functions $\varphi:[1,\infty)\rightarrow(0,\infty)$ such that both the functions $\varphi$ and $1/\varphi$ are bounded on each compact interval $[1,b]$, with $1<b<\infty$, and that the function $\varphi$ is slowly varying at infinity in the sense of J.~Karamata~\cite{Karamata30a}. The latter property means that $\varphi(\lambda t)/\varphi(t)\rightarrow 1$ as $t\rightarrow\infty$ for every  $\lambda>0$.

Note that the slowly varying functions are well investigated and has various applications (see monographs \cite{Seneta76, BinghamGoldieTeugels89}). We obtain an important example of these functions by setting
$$
\varphi(t):=(\log t)^{r_{1}}(\log\log
t)^{r_{2}}\ldots(\underbrace{\log\ldots\log}_{k\;\mbox{\small times}}
t)^{r_{k}}\quad\mbox{of}\quad t\gg1.
$$
Here, $k\in\mathbb{Z}$, with $k\geq1$, and $r_{1},\ldots,r_{k}\in\mathbb{R}$ are arbitrary parameters.

Let $s\in\mathbb{R}$ and $\varphi\in\mathcal {M}$. By definition, the complex linear space $H^{s,\varphi}(\mathbb{R}^{n})$, with $n\geq1$, consists of all distributions $w\in\mathcal{S}'(\mathbb{R}^{n})$ that the Fourier transform $\widehat{w}$ of $w$ is locally Lebesgue integrable over $\mathbb{R}^{n}$ and satisfies the condition
$$
\int\limits_{\mathbb{R}^{n}}
\langle\xi\rangle^{2s}\varphi^{2}(\langle\xi\rangle)\,
|\widehat{w}(\xi)|^{2}\,d\xi<\infty.
$$
Here, as usual, $\mathcal{S}'(\mathbb{R}^{n})$ is the linear topological space of all tempered distributions on $\mathbb{R}^{n}$, and
$\langle\xi\rangle:=(1+|\xi|^{2})^{1/2}$ is the smoothed modulus of $\xi\in\mathbb{R}^{n}$. The space $H^{s,\varphi}(\mathbb {R}^{n})$ is endowed with the inner product
$$
(w_{1},w_{2})_{H^{s,\varphi}(\mathbb {R}^{n})}:=
\int\limits_{\mathbb{R}^{n}}
\langle\xi\rangle^{2s}\,\varphi^{2}(\langle\xi\rangle)\,
\widehat{w_{1}}(\xi)\,\overline{\widehat{w_{2}}(\xi)}\,d\xi.
$$
It naturally induces the norm
$$\|w\|_{H^{s,\varphi}(\mathbb{R}^{n})}:=
(w,w)_{H^{s,\varphi}(\mathbb{R}^{n})}^{1/2}.
$$

The space $H^{s,\varphi}(\mathbb{R}^{n})$ just defined is a special isotropic Hilbert case of the spaces $\mathcal{B}_{p,\mu}$ introduced and investigated by L.~H\"ormander in \cite{Hermander63} (Sect.~2.2) (see also his monograph \cite{Hermander83} (Sect.~10.1)). Namely, $H^{s,\varphi}(\mathbb{R}^{n})=\mathcal{B}_{p,\mu}$ if $p=2$ and $\mu(\xi)\equiv\langle\xi\rangle^{s}\varphi(\langle\xi\rangle)$. Note that the inner product spaces $\mathcal{B}_{2,\mu}$ were also investigated by L.~R.~Volevich and B.~P.~Paneah in \cite{VolevichPaneah65} (\S~2).

In the $\varphi\equiv1$ case, the space $H^{s,\varphi}(\mathbb {R}^{n})$ becomes the inner product Sobolev space $H^{s}(\mathbb {R}^{n})$ of order $s\in\mathbb{R}$. Generally, we have the continuous and dense embeddings
\begin{equation}\label{3f7}
H^{s+\varepsilon}(\mathbb{R}^{n})\hookrightarrow H^{s,\varphi}(\mathbb{R}^{n})\hookrightarrow H^{s-\varepsilon}(\mathbb{R}^{n})
\quad\mbox{for every}\quad\varepsilon>0.
\end{equation}
They show that, the numerical parameter $s$ characterizes the main regularity of the distributions $w\in H^{s,\varphi}(\mathbb{R}^{n})$, whereas the function parameter $\varphi$ defines a supplementary regularity, which is subordinate to $s$. Specifically, if $\varphi(t)\rightarrow\infty$ [or $\varphi(t)\rightarrow 0$] as $t\rightarrow\infty$, then the parameter $\varphi$ defines the supplementary positive [or negative] regularity. Thus, we can say that $\varphi$ refines the main regularity $s$ in the class of Hilbert separable spaces
$$
\bigl\{H^{s,\varphi}(\mathbb{R}^{n}):\,
s\in\mathbb{R},\,\varphi\in\mathcal{M}\bigr\}.
$$
This class is selected by V.~A.~Mikhailets and A.~A.~Murach in \cite{MikhailetsMurach05UMJ5} and is called the refined Sobolev scale over $\mathbb{R}^{n}$ (see also monograph \cite{MikhailetsMurach14} (Subsect.~1.3.3)).

In the paper, we need the refined Sobolev scales over the domain $\Omega$ and its boundary~$\Gamma$. These scales are constructed in the standard way with the help of the H\"ormander spaces $H^{s,\varphi}(\mathbb{R}^{n})$. Let us pass to the corresponding definitions \cite{MikhailetsMurach14} (Sect.~2.1 and 3.1).

The complex linear space $H^{s,\varphi}(\Omega)$ and the norm in it are defined as follows:
\begin{gather*}
H^{s,\varphi}(\Omega):=\bigl\{w\!\upharpoonright\!\Omega:\,
w\in H^{s,\varphi}(\mathbb{R}^{n})\bigr\},\\
\|u\|_{H^{s,\varphi}(\Omega)}:=
\inf\,\bigl\{\|w\|_{H^{s,\varphi}(\mathbb{R}^{n})}:\,
w\in H^{s,\varphi}(\mathbb{R}^{n}),\;\,
u=w\!\upharpoonright\!\Omega\bigr\},
\end{gather*}
with $u\in H^{s,\varphi}(\Omega)$. Here, as usual, $w\!\upharpoonright\!\Omega$ stands for the restriction of the distribution $w$ to the domain $\Omega$. The space $H^{s,\varphi}(\Omega)$ is Hilbert and separable with respect to this norm. The space $H^{s,\varphi}(\Omega)$ is continuously embedded in the linear topological space $\mathcal{S}'(\Omega):=
\{w\!\upharpoonright\!\Omega:w\in\mathcal{S}'(\mathbb{R}^{n})\}$.
The set $C^{\infty}(\overline{\Omega})$ is dense in $H^{s,\varphi}(\Omega)$.

In short, the space $H^{s,\varphi}(\Gamma)$ consists of all distributions on $\Gamma$ that belong, in local coordinates, to $H^{s,\varphi}(\mathbb{R}^{n-1})$. Let us dwell on the definition of $H^{s,\varphi}(\Gamma)$. From $C^{\infty}$-structure on $\Gamma$, we arbitrarily choose a finite collection of local charts $\alpha_j: \mathbb{R}^{n-1}\leftrightarrow\Gamma_{j}$, with $j=1,\ldots,\lambda$, that the open sets $\Gamma_{1},\ldots,\Gamma_{\lambda}$ form a covering of the manifold $\Gamma$. Besides, we arbitrarily choose functions $\chi_j\in C^{\infty}(\Gamma)$, with $j=1,\ldots,\lambda$, that form a partition of unity on $\Gamma$ which satisfies the condition $\mathrm{supp}\,\chi_j\subset\Gamma_j$. Then the complex linear space $H^{s,\varphi}(\Gamma)$ and the norm in it are defined as follows:
\begin{gather*}
H^{s,\varphi}(\Gamma):=\bigl\{h\in\mathcal{D}'(\Gamma):\,
(\chi_{j}h)\circ\alpha_{j}\in H^{s,\varphi}(\mathbb{R}^{n-1})\;\,\mbox{for every}\;\,j\in\{1,\ldots,\lambda\}\bigr\},\\
\|h\|_{H^{s,\varphi}(\Gamma)}:=\biggl(\;\sum_{j=1}^{\lambda}\,
\|(\chi_{j}h)\circ\alpha_{j}\|_{H^{s,\varphi}(\mathbb{R}^{n-1})}^{2}
\biggr)^{1/2}.
\end{gather*}
Here, as usual, $\mathcal{D}'(\Gamma)$ is the linear topological space of all distributions on $\Gamma$, whereas $(\chi_{j}h)\circ\alpha_{j}$ is the representation of the distribution $\chi_{j}h$ in the local chart $\alpha_{j}$. The space $H^{s,\varphi}(\Gamma)$ is Hilbert and separable with respect to this norm. This space does not depend (up to equivalence of norms) on the indicated choice of the local charts and partition of unity \cite{MikhailetsMurach14} (Theorem~2.3). The space $H^{s,\varphi}(\Gamma)$ is continuously embedded in $\mathcal{D}'(\Gamma)$. The set $C^{\infty}(\Gamma)$ is dense in $H^{s,\varphi}(\Gamma)$.

Thus, we have the refined Sobolev scales
\begin{equation}\label{3f8}
\bigl\{H^{s,\varphi}(\Omega):\,
s\in\mathbb{R},\,\varphi\in\mathcal{M}\bigr\}\quad\mbox{and}\quad
\bigl\{H^{s,\varphi}(\Gamma):\,
s\in\mathbb{R},\,\varphi\in\mathcal{M}\bigr\}
\end{equation}
over $\Omega$ and $\Gamma$ respectively. They contain the inner product Sobolev spaces $H^{s}(\Omega):=H^{s,1}(\Omega)$ and $H^{s}(\Gamma):=H^{s,1}(\Gamma)$. (Here, of course, $1$ denotes the function equaled identically to $1$). It follows from \eqref{3f7} that
\begin{gather}\label{3f9}
H^{s+\varepsilon}(\Omega)\hookrightarrow H^{s,\varphi}(\Omega)\hookrightarrow H^{s-\varepsilon}(\Omega)
\quad\mbox{for every}\quad\varepsilon>0,\\
H^{s+\varepsilon}(\Gamma)\hookrightarrow H^{s,\varphi}(\Gamma)\hookrightarrow H^{s-\varepsilon}(\Gamma)
\quad\mbox{for every}\quad\varepsilon>0.
\label{3f10}
\end{gather}
These embeddings are compact and dense.

It is worthwhile to note the following connection between the scales \eqref{3f8} (see \cite{MikhailetsMurach14} (Subsect. 3.2.1)). Let $s>1/2$ and $\varphi\in\mathcal{M}$, then the trace mapping $u\mapsto u\!\upharpoonright\!\Gamma$, with $u\in C^{\infty}(\Gamma)$, extends uniquely (by continuity) to a bounded surjective operator $R_{\Gamma}:H^{s,\varphi}(\Omega)\rightarrow H^{s-1/2,\varphi}(\Gamma)$. Thus, for each distribution $u\in H^{s,\varphi}(\Omega)$, its trace $R_{\Gamma}u$ on $\Gamma$ is well defined. Moreover, we have the equivalence of norms
$$
\|h\|_{H^{s-1/2,\varphi}(\Gamma)}\asymp
\inf\,\bigl\{\|u\|_{H^{s,\varphi}(\Omega)}:\,
u\in H^{s,\varphi}(\Omega),\;\,h=R_{\Gamma}u\bigr\}
$$
on the class of all $h\in H^{s-1/2,\varphi}(\Gamma)$. But it is impossible to define this trace reasonably in the $s<1/2$ case.

\section{The main results} 

To formulate them, we need to introduce the linear spaces $N$ and $N^{+}$ of infinitely smooth solutions to the problem \eqref{3f1}, \eqref{3f2} and its formally adjoint problem \eqref{3f4}, \eqref{3f5},~\eqref{3f6}. Namely, let $N$ consist of all solutions
$$
(u,v_{1},...,v_{\varkappa})\in C^{\infty}(\overline{\Omega})\times
\bigl(C^{\infty}(\Gamma)\bigr)^{\varkappa}
$$
to the problem \eqref{3f1}, \eqref{3f2} in the case where $f=0$ on $\Omega$ and all $g_{j}=0$ on~$\Gamma$. Besides, let $N^{+}$ consist of all solutions
$$
(w,h_{1},...,h_{q+\varkappa})\in C^{\infty}(\overline{\Omega})\times
\bigl(C^{\infty}(\Gamma)\bigr)^{q+\varkappa}
$$
to the problem \eqref{3f4}, \eqref{3f5}, \eqref{3f6} in the case where $\omega=0$ on $\Omega$ and all $\chi_{j}=0$ and all $\chi_{2q+j}=0$ on $\Gamma$. Since these problems are elliptic on $\Omega$, the spaces $N$ and $N^{+}$ are finite-dimensional \cite{KozlovMazyaRossmann97} (Lemma~3.4.2).

We are motivated by the following result.

\textbf{Proposition 1.} \it For arbitrary $s>2q-1/2$ and $\varphi\in\mathcal{M}$, the mapping \eqref{3f3} extends uniquely (by continuity) to a bounded operator
\begin{equation}\label{3f11}
\Lambda:\,H^{s,\varphi}(\Omega)\oplus
\bigoplus_{k=1}^{\varkappa}H^{s+r_{k}-1/2,\varphi}(\Gamma)\rightarrow
H^{s-2q,\varphi}(\Omega)\oplus
\bigoplus_{j=1}^{q+\varkappa}H^{s-m_{j}-1/2,\varphi}(\Gamma).
\end{equation}
This operator is Fredholm. Its kernel coincides with $N$, and its range consists of all vectors
$$
(f,g_{1},\ldots,g_{q+\varkappa})\in
H^{s-2q,\varphi}(\Omega)\oplus
\bigoplus_{j=1}^{q+\varkappa}H^{s-m_{j}-1/2,\varphi}(\Gamma)
$$
such that
\begin{equation}\label{3f12}
(f,w)_{\Omega}+\sum_{j=1}^{q+\varkappa}(g_{j},h_{j})_{\Gamma}=0\quad
\mbox{for all}\quad(w,h_{1},\ldots,h_{q+\varkappa})\in N^{+}.
\end{equation}
The index of the operator \eqref{3f11} is equal to $\dim N-\dim N^{+}$ and does not depend on $s$ and~$\varphi$. \rm

Here and below, $(\cdot,\cdot)_{\Omega}$ and $(\cdot,\cdot)_{\Gamma}$ denote the extension by continuity of the inner products in the Hilbert spaces $L_{2}(\Omega)$ and $L_{2}(\Gamma)$ respectively.

In the Sobolev case of $\varphi(t)\equiv1$, Proposition~1 is proved in \cite{KozlovMazyaRossmann97} (Theorem~3.4.1) for integer $s$ and in  \cite{Roitberg99} (Theorems 2.4.1 and 2.7.3) for real $s$ and general elliptic systems. For arbitrary $\varphi\in\mathcal{M}$ and real $s>2q$, this proposition is proved in \cite{Chepurukhina14Coll2} (Theorem~1). The case where $\varphi\in\mathcal{M}$ and $2q-1/2<s\leq2q$ is covered in \cite{ChepurukhinaMurach14MFAT1} (Theorem~1).

In connection with Proposition~1, we recall the definition of a Fredholm operator and its index. A linear bounded operator $T:E_{1}\rightarrow E_{2}$ acting between certain Banach spaces $E_{1}$ and $E_{2}$ is said to be Fredholm if its kernel $\ker T$ and co-kernel $E_{2}/T(E_{1})$ are finite-dimensional. If this operator is Fredholm, then its range $T(E_{1})$ is closed in $E_{2}$; see, e.g., \cite{Hermander85} (Lemma~19.1.1). The number  $\mathrm{ind}\,T:=\dim\ker T-\dim(E_{2}/T(E_{1}))$ is called the index of the Fredholm operator~$T$.

Proposition~1 is not true if $s\leq2q-1/2$. This is caused by the fact that the boundary differential operators $B_{j}$ cannot be reasonably applied to all distributions $u\in H^{s,\varphi}(\Omega)$ in this case. Namely, the mapping $u\mapsto B_{j}u$, with $u\in C^{\infty}(\overline{\Omega})$, cannot be extended to a continuous operator from the whole $H^{s,\varphi}(\Omega)$ to $\mathcal{D}'(\Gamma)$ if $s\leq m_{j}+1/2$.

To obtain a version of Theorem~1 for arbitrary $s\leq2q-1/2$, we restrict ourselves to the solutions $u\in H^{s,\varphi}(\Omega)$ of the elliptic equation $Au=f$ with $f\in L_{2}(\Omega)=H^{0}(\Omega)$. This approach originates from J.-L.~Lions and E.~Magenes \cite{LionsMagenes62, LionsMagenes63}, who took it to the investigation of regular elliptic boundary-value problems in Sobolev spaces.

Let $s<2q$ and $\varphi\in\mathcal{M}$. Consider the linear space
$$
H^{s,\varphi}_{A}(\Omega):=\bigl\{u\in H^{s,\varphi}(\Omega):\,
Au\in L_{2}(\Omega)\bigr\}
$$
endowed with the graph norm
\begin{equation}\label{3f13}
\|u\|_{H^{s,\varphi}_{A}(\Omega)}:=
\bigl(\,\|u\|^{2}_{H^{s,\varphi}(\Omega)}+
\|Au\|^{2}_{L_{2}(\Omega)}\bigr)^{1/2}.
\end{equation}
Here, $Au$ is understood in the sense of the theory of distributions.

This space is Hilbert with respect to the norm \eqref{3f13}. Indeed,
this norm is generated by an inner product because so are the norms on the right of \eqref{3f13}. Besides, the space $H^{s,\varphi}_{A}(\Omega)$ is complete with respect to this norm. In fact, if $(u_{k})$ is a Cauchy sequence in this space, then there exist the limits $u:=\lim u_{k}$ in $H^{s,\varphi}(\Omega)$ and $f:=\lim Au_{k}$ in $L_{2}(\Omega)$
because the spaces $H^{s,\varphi}(\Omega)$ and $L_{2}(\Omega)$ are complete. The differential operator $A$ is continuous in $\mathcal{S}'(\Omega)$, so
$Au=\lim Au_{k}=f$ in $\mathcal{S}'(\Omega)$. Here, recall, $u\in H^{s,\varphi}(\Omega)$ and $f\in L_{2}(\Omega)$. Therefore, $u\in H^{s,\varphi}_{A}(\Omega)$ and $\lim u_{k}=u$ in the space $H^{s,\varphi}_{A}(\Omega)$. Thus, this space is complete.

It is useful to note that even when all coefficients of $A$ are constant, the space $H^{s,\varphi}_{A}(\Omega)$ depends essentially on each of them. This follows from the result by L.~H\"ormander \cite{Hermander55} (Theorem~3.1).

Consider the Hilbert spaces
\begin{gather*}
\mathcal{D}^{s,\varphi}_{A}(\Omega,\Gamma):=H^{s,\varphi}_{A}(\Omega)\oplus
\bigoplus_{k=1}^{\varkappa}H^{s+r_{k}-1/2,\varphi}(\Gamma),\\
\mathcal{E}_{0,s,\varphi}(\Omega,\Gamma):=L_{2}(\Omega)\oplus
\bigoplus_{j=1}^{q+\varkappa}H^{s-m_{j}-1/2,\varphi}(\Gamma).
\end{gather*}
We investigate the properties of the elliptic boundary-value problem \eqref{3f1}, \eqref{3f2} in the couple of these spaces.

\textbf{Theorem 1.} \it Let $s<2q$ and $\varphi\in\mathcal{M}$. Then the set $C^{\infty}(\overline{\Omega})$ is dense in the space $H^{s,\varphi}_{A}(\Omega)$, and the mapping \eqref{3f3} extends uniquely (by continuity) to a bounded operator
\begin{equation}\label{3f14}
\Lambda:\,\mathcal{D}^{s,\varphi}_{A}(\Omega,\Gamma)\rightarrow
\mathcal{E}_{0,s,\varphi}(\Omega,\Gamma).
\end{equation}
This operator is Fredholm. Its kernel coincides with $N$, whereas its range consists of all vectors
\begin{equation}\label{3f15}
(f,g):=(f,g_{1},\ldots,g_{q+\varkappa})\in
\mathcal{E}_{0,s,\varphi}(\Omega,\Gamma)
\end{equation}
that satisfy condition \eqref{3f12}. The index of the operator \eqref{3f14} is equal to $\dim N-\dim N^{+}$ and does not depend on $s$ and~$\varphi$. \rm

If $N=\{0\}$ and $N^{+}=\{0\}$, then the operator \eqref{3f14} is an isomorphism between the spaces $\mathcal{D}^{s,\varphi}_{A}(\Omega,\Gamma)$ and $\mathcal{E}_{0,s,\varphi}(\Omega,\Gamma)$. Generally, this operator naturally induces an isomorphism
\begin{equation}\label{3f-isom}
\Lambda:\,\mathcal{D}^{s,\varphi}_{A}(\Omega,\Gamma)/N\leftrightarrow
\mathcal{R}_{0,s,\varphi}(\Omega,\Gamma).
\end{equation}
Here,
\begin{equation*}
\mathcal{R}_{0,s,\varphi}(\Omega,\Gamma):=
\bigl\{(f,g_{1},\ldots,g_{q+\varkappa})\in
\mathcal{E}_{0,s,\varphi}(\Omega,\Gamma):\,\mbox{\eqref{3f12} is true}\bigr\}
\end{equation*}
is a subspace of $\mathcal{E}_{0,s,\varphi}(\Omega,\Gamma)$.

Let us discuss the properties of generalized solutions to the elliptic boundary-value problem \eqref{3f1}, \eqref{3f2}. We first give the definition of such solutions.

Let us denote
$$
\mathcal{S}'_{A}(\Omega):=\{u\in\mathcal{S}'(\Omega):\,
Au\in L_{2}(\Omega)\}.
$$
We interpret $\mathcal{S}'_{A}(\Omega)$ as a subspace of the linear topological space $\mathcal{S}'(\Omega)$. Since the domain $\Omega$ is bounded, we conclude that $\mathcal{S}'_{A}(\Omega)$ is the union of all spaces $H^{s,\varphi}_{A}(\Omega)$ with $s<2q$ and $\varphi\in\mathcal{M}$.

Consider an arbitrary vector
\begin{equation}\label{3f19}
(u,v):=(u,v_{1},\ldots,v_{\varkappa})\in
\mathcal{S}'_{A}(\Omega)\times(\mathcal{D}'(\Gamma))^{\varkappa}.
\end{equation}
Since $(u,v)\in\mathcal{D}^{s,\varphi}_{A}(\Omega,\Gamma)$ for certain  $s<2q$ and $\varphi\in\mathcal{M}$, we define its image
\eqref{3f15} by the formula $\Lambda(u,v)=(f,g)$ if we apply the bounded operator \eqref{3f14}. The vector $(u,v)$ does not depend on $s$ and $\varphi$ and is called a generalized solution of the boundary-value problem \eqref{3f1}, \eqref{3f2} with the right-hand sides \eqref{3f15}. This solution satisfies the following a priori estimate.

\textbf{Theorem 2.} \it Let $s<2q$, $\varphi\in\mathcal{M}$, and  $\sigma>0$. Then there exists a number $c>0$ such that
\begin{equation}\label{3f20}
\|(u,v)\|_{\mathcal{D}^{s,\varphi}_{A}(\Omega,\Gamma)}\leq c\,\bigl(\,\|\Lambda(u,v)\|_{\mathcal{E}_{0,s,\varphi}(\Omega,\Gamma)}+
\|(u,v)\|_{\mathcal{D}_A^{s-\sigma,\varphi}(\Omega,\Gamma)}\bigl)
\end{equation}
for an arbitrary vector  $(u,v)\in\mathcal{D}^{s,\varphi}_{A}(\Omega,\Gamma)$. Here,  $c=c(s,\varphi,\sigma)$ is independent of $(u,v)$. \rm

If $N=\{0\}$, then the second term on the right of the inequality \eqref{3f20} is absent. This immediately follows from the isomorphism \eqref{3f-isom}.

In connection with this theorem, we mention Peetre's lemma \cite{Peetre61} (Lemma~3) and remark that the embedding $\mathcal{D}^{s,\varphi}_{A}(\Omega,\Gamma)\hookrightarrow\mathcal{D}_A^{s-\sigma,\varphi}(\Omega,\Gamma)$, with $\sigma>0$, is not compact.

The next theorem says about the regularity properties of generalized solutions to the elliptic boundary-value problem \eqref{3f1}, \eqref{3f2}.

\textbf{Theorem 3.} \it Let $s\in\mathbb{R}$ and $\varphi\in\mathcal{M}$. Suppose that the vector \eqref{3f19} is a generalized solution to the elliptic boundary-value problem \eqref{3f1}, \eqref{3f2} whose right-hand sides meets the conditions
\begin{gather}\label{3f21}
f\in H^{s-2q,\varphi}(\Omega)
\quad\mbox{if}\quad s\geq2q,\\
g_{j}\in H^{s-m_{j}-1/2,\varphi}(\Gamma)
\quad\mbox{for each}\quad j\in\{1,...,q+\varkappa\}. \label{3f22}
\end{gather}
Then the solution satisfies the conditions
\begin{gather}\label{3f-a}
u\in H^{s,\varphi}(\Omega),\\
v_{k}\in H^{s+r_{k}-1/2,\varphi}(\Gamma)
\quad\mbox{for each}\quad k\in\{1,...,\varkappa\}.\label{3f-b}
\end{gather}
\rm

We conclude that the solution inherits the (supplementary) regularity $\varphi$ of the right-hand sides of the elliptic problem.

We note that the restriction $s\geq2q$ in condition \eqref{3f21} is caused by the following: if $s<2q$, then
$$
f=Au\in L_{2}(\Omega)\subset H^{s-2q,\varphi}(\Omega)
$$
in view of \eqref{3f19} and \eqref{3f9}. Thus, the assumption $f\in H^{s-2q,\varphi}(\Omega)$ is superfluous in Theorem~3 in the case of $s<2q$.

\section{Auxiliary results relating to the problem under consideration} 

Our proof of Theorem~1 will be based on the result by V.~A.~Kozlov, V.~G.~Maz'ya and J.~Rossmann \cite{KozlovMazyaRossmann97} (Theorem~3.4.1). They demonstrated that the elliptic boundary-value problem \eqref{3f1}, \eqref{3f2} is Fredholm in the two-sided scale of function spaces introduced by Ya.~A.~Roitberg \cite{Roitberg64, Roitberg65}. Being a certain modification of the classical Sobolev scale, this two-sided scale proves to be useful in the theory of elliptic boundary-value problems in the situation where their right-hand sides can be nonregular distributions (see monographs by Yu.~M.~Berezansky \cite{Berezansky68} (Chapt.~III, Sect.~6) and Ya.~A.~Roitberg \cite{Roitberg96, Roitberg99}).

Let us recall the definition of the inner product spaces $H^{s,(2q)}(\Omega)$ introduced by Ya.~A.~Roitberg. For our purposes, we restrict ourselves to the case of $s\in\mathbb{Z}$.

Beforehand, we need to define the space $H^{s,(0)}(\Omega)$. If $s\geq0$, then $H^{s,(0)}(\Omega):=H^{s}(\Omega)$ is the inner product Sobolev space of order $s$ in $\Omega$. But, if $s<0$, then $H^{s,(0)}(\Omega)$ is the completion of $C^{\infty}(\overline{\Omega})$ with respect to the Hilbert norm
$$
\|u\|_{H^{s,(0)}(\Omega)}:=
\sup\left\{\,\frac{|(u,w)_{\Omega}|}
{\;\quad\quad\quad\|w\|_{H^{-s}(\Omega)}}\,:
\,w\in H^{-s}(\Omega),\,w\neq0\right\}.
$$
Thus, we have the rigging
$$
H^{s,(0)}(\Omega)\hookrightarrow L_{2}(\Omega)\hookrightarrow H^{-s,(0)}(\Omega),\quad\mbox{with}\quad s>0,
$$
of the space $L_{2}(\Omega)$ with positive and negative Sobolev spaces (see \cite{Berezansky68} (Chapt.~I, Sect.~3, Subsect.~1)). Note that the negative space $H^{s,(0)}(\Omega)$, with $s<0$, admits the following description: the mapping $u\mapsto\mathcal{O}u$, with $u\in C^{\infty}(\overline{\Omega})$, extends uniquely (by continuity) to an isometric isomorphism of the space $H^{s,(0)}(\Omega)$ onto the subspace $\{w\in H^{s}(\mathbb{R}^{n}):\mathrm{supp}\,w\subseteq \overline{\Omega}\,\}$ of $H^{s}(\mathbb{R}^{n})$. Here, $\mathcal{O}u:=u$ on $\overline{\Omega}$, and $\mathcal{O}u:=0$ on $\mathbb{R}^{n}\setminus\overline{\Omega}$.

We can now define the Hilbert space $H^{s,(2q)}(\Omega)$, with $s\in\mathbb{Z}$. By definition, the space $H^{s,(2q)}(\Omega)$ is the  completion of the set $C^{\infty}(\overline{\Omega})$ with respect to the Hilbert norm
$$
\|u\|_{H^{s,(2q)}(\Omega)}:=
\biggl(\|u\|_{H^{s,(0)}(\Omega)}^{2}+
\sum_{k=1}^{2q}\;\|(D_{\nu}^{k-1}u)\!\upharpoonright\!\Gamma\|
_{H^{s-k+1/2}(\Gamma)}^{2}\biggr)^{1/2}.
$$

This space admits the following description \cite{Roitberg96} (Sect. 2.2): the linear mapping
$$
T_{2q}:u\mapsto\bigl(u,u\!\upharpoonright\!\Gamma,\ldots,
(D_{\nu}^{2q-1}u)\!\upharpoonright\!\Gamma\bigr),\quad\mbox{with}\quad u\in C^{\infty}(\overline{\Omega}),
$$
extends uniquely (by continuity) to an isometric operator
$$
T_{2q}:\,H^{s,(2q)}(\Omega)\rightarrow
H^{s,(0)}(\Omega)\oplus
\bigoplus_{k=1}^{2q}\,H^{s-k+1/2}(\Gamma)=:
\Pi_{s,(2q)}(\Omega,\Gamma);
$$
the range of this operator consists of all vectors
$$
(u_{0},u_{1},\ldots,u_{2q})\in\Pi_{s,(2q)}(\Omega,\Gamma)
$$
such that $u_{k}=R_{\Gamma}D_{\nu}^{k-1}u_{0}$ for all integers $k\in\{1,\ldots,2q\}$ satisfying $s>k-1/2$. Thus, we can interpret an arbitrary element $u\in H^{s,(2q)}(\Omega)$ as the vector $(u_{0},u_{1},\ldots,u_{2q}):=T_{2q}u$. This interpretation is taken in the mentioned monograph \cite{KozlovMazyaRossmann97} (Sect. 3.2.1) as a definition of the space $H^{s,(2q)}(\Omega)$.

Note that
\begin{equation}\label{3f25}
H^{s,(2q)}(\Omega)=H^{s}(\Omega)\quad\mbox{for all}\quad s\geq 2q
\end{equation}
as completions of $C^{\infty}(\overline{\Omega})$ with respect to equivalent norms.

We can now formulate the result by V.~A.~Kozlov, V.~G.~Maz'ya and J.~Rossmann \cite{KozlovMazyaRossmann97} (Theorem~3.4.1). Consider the Hilbert spaces
\begin{gather*}
\mathcal{D}^{s,(2q)}(\Omega,\Gamma):=H^{s,(2q)}(\Omega)\oplus
\bigoplus_{k=1}^{\varkappa}H^{s+r_{k}-1/2}(\Gamma),\\
\mathcal{E}^{s-2q,(0)}(\Omega,\Gamma):=H^{s-2q,(0)}(\Omega)\oplus
\bigoplus_{j=1}^{q+\varkappa}H^{s-m_{j}-1/2}(\Gamma),
\end{gather*}
with $s\in\mathbb{Z}$.

\textbf{Proposition 2.} \it For arbitrary $s\in\mathbb{Z}$, the mapping \eqref{3f3} extends uniquely (by continuity) to a bounded operator
\begin{equation}\label{3f26}
\Lambda:\,\mathcal{D}^{s,(2q)}(\Omega,\Gamma)\rightarrow
\mathcal{E}^{s-2q,(0)}(\Omega,\Gamma).
\end{equation}
This operator is Fredholm. Its kernel coincides with $N$, and its range consists of all vectors
\begin{equation}\label{3f27}
(f,g_{1},\ldots,g_{q+\varkappa})\in\mathcal{E}^{s-2q,(0)}(\Omega,\Gamma)
\end{equation}
that satisfy \eqref{3f12}. The index of the operator \eqref{3f26} is equal to $\dim N-\dim N^{+}$ and does not depend on $s$ and~$\varphi$. \rm

Note that if $s\in\mathbb{Z}$, $s\geq2q$, and $\varphi\equiv1$, then Proposition~1 coincides with Proposition~2 in view of \eqref{3f25}.

In the proof of Theorem~1, we will use Proposition~2 together with the following result. For an arbitrary integer $s<2q$, we consider the linear space
$$
H^{s,(2q)}_{A}(\Omega):=\bigl\{u\in H^{s,(2q)}(\Omega):\,
Au\in L_{2}(\Omega)\bigr\}
$$
endowed with the graph norm
\begin{equation}\label{3f28}
\|u\|_{H^{s,(2q)}_{A}(\Omega)}:=
\bigl(\,\|u\|^{2}_{H^{s,(2q)}(\Omega)}+
\|Au\|^{2}_{L_{2}(\Omega)}\bigr)^{1/2}.
\end{equation}
Note that the image $Au\in H^{s-2q,(0)}(\Omega)$ is well defined for each
$u\in H^{s,(2q)}(\Omega)$ by closure because the mapping $u\mapsto Au$, with $u\in C^{\infty}(\overline{\Omega})$, extends by continuity to a bounded operator $A:H^{s,(2q)}(\Omega)\rightarrow H^{s-2q,(0)}(\Omega)$ for every $s\in\mathbb{Z}$ \cite{Roitberg96} (Lemma 2.3.1). The space $H^{s,(2q)}_{A}(\Omega)$ is Hilbert with respect to the norm \eqref{3f28}, and the set $C^{\infty}(\overline{\Omega})$ is dense in it (see, e.g., \cite{MikhailetsMurach14} (Sect. 4.4.2, the proof of Theorem 4.25)).

\textbf{Proposition 3.} \it For an arbitrary integer $s<2q$, the identity mapping on $C^{\infty}(\overline{\Omega})$ extends uniquely (by continuity) to an isomorphism between the spaces $H^{s,(2q)}_{A}(\Omega)$ and $H^{s}_{A}(\Omega)$. Thus, these spaces are equal as completions of the set $C^{\infty}(\overline{\Omega})$ with respect to equivalent norms. \rm

Here and below, $H^{s}_{A}(\Omega)$ is the space $H^{s,\varphi}_{A}(\Omega)$ in the Sobolev case of $\varphi\equiv1$.

The proof of Proposition 3 is given in monograph \cite{MikhailetsMurach14} (Sect. 4.4.2, see the reasoning from formula (4.196) to the end of the proof of Theorem~4.25).

\section{Auxiliary results relating to\\interpolation between Hilbert spaces} 

The refined Sobolev scales possess an important interpolation property with respect to inner product Sobolev spaces. Namely, each space $H^{s,\varphi}(G)$, with $s\in\mathbb{R}$, $\varphi\in\mathcal{M}$, and $G\in\{\mathbb{R}^{n},\Omega,\Gamma\}$, can be obtained by the interpolation with an appropriate function parameter between the Sobolev spaces $H^{s-\varepsilon}(G)$ and $H^{s+\delta}(G)$ with $\varepsilon,\delta>0$. This property will play a key role in our proof of Theorem~1 for arbitrary $\varphi\in\mathcal{M}$. Therefore, we recall the definition of this interpolation in the case of arbitrary Hilbert spaces and discuss some of its properties that are necessary for us.

The method of interpolation with a function parameter between Hilbert spaces appeared first in C.~Foia\c{s} and J.-L.~Lions' paper \cite[p.~278]{FoiasLions61}. We follow monograph \cite{MikhailetsMurach14} (Sect.~1.1), which systematically expound this method. For our purposes, it is sufficient to restrict ourselves to separable Hilbert spaces.

Let $X:=[X_{0},X_{1}]$ be a given ordered couple of separable complex Hilbert spaces $X_{0}$ and $X_{1}$ such that the continuous and dense embedding $X_{1}\hookrightarrow X_{0}$ holds. This couple is said to be admissible. For $X$, there exists a self-adjoint positive-definite operator $J$ on $X_{0}$ that has the domain $X_{1}$ and that satisfies the condition $\|Jw\|_{X_{0}}=\|w\|_{X_{1}}$ for each $w\in X_{1}$. This operator is uniquely determined by the couple $X$ and is called a generating operator for~$X$. It defines an isometric isomorphism $J:X_{1}\leftrightarrow X_{0}$.

Let $\mathcal{B}$ denote the set of all Borel measurable functions
$\psi:(0,\infty)\rightarrow(0,\infty)$ such that $\psi$ is bounded on each compact interval $[a,b]$, with $0<a<b<\infty$, and that $1/\psi$ is bounded on every set $[r,\infty)$, with $r>0$.

For arbitrary $\psi\in\mathcal{B}$, we consider the (generally, unbounded) operator $\psi(J)$, which is defined on $X_{0}$ as the Borel function $\psi$ of $J$ and is built with the help of Spectral Theorem applied to the self-adjoint operator $J$. Let $[X_{0},X_{1}]_{\psi}$ or, simply, $X_{\psi}$ denote the domain of the operator $\psi(J)$ endowed with the norm $\|u\|_{X_{\psi}}=\|\psi(J)u\|_{X_{0}}$. The space $X_{\psi}$ is Hilbert and separable with respect to this norm.

A function $\psi\in\mathcal{B}$ is called an interpolation parameter if the following condition is fulfilled for all admissible couples $X=[X_{0},X_{1}]$ and $Y=[Y_{0},Y_{1}]$ of Hilbert spaces and for an arbitrary linear mapping $T$ given on $X_{0}$: if the restriction of $T$ to $X_{j}$ is a bounded operator $T:X_{j}\rightarrow Y_{j}$ for each $j\in\{0,1\}$, then the restriction of $T$ to $X_{\psi}$ is also a bounded operator $T:X_{\psi}\rightarrow Y_{\psi}$. In this case we say that $X_{\psi}$ is obtained by the interpolation with the function parameter $\psi$ of the couple $X$ (or between the spaces $X_{0}$ and $X_{1}$).

The function $\psi\in\mathcal{B}$ is an interpolation parameter if and only if $\psi$ is pseudoconcave on a neighbourhood of infinity, i.e.
$\psi(t)\asymp\psi_{1}(t)$ with $t\gg1$ for a certain positive concave
function~$\psi_{1}(t)$. (As usual, $\psi\asymp\psi_{1}$ means that both the functions $\psi/\psi_{1}$ and $\psi_{1}/\psi$ are bounded on the indicated set). This fundamental property of the interpolation follows from J.~Peetre's results \cite{Peetre66, Peetre68}. Specifically, if
$\psi(t)\equiv t^{\sigma}\psi_{0}(t)$ for some number $\sigma\in(0,1)$ and function $\psi_{0}$ varying slowly at infinity (in the sense of J.~Karamata), then $\psi$ is an interpolation parameter.

Let us formulate the above-mentioned interpolation property of the refined Sobolev scales \cite{MikhailetsMurach14} (Theorems 1.14, 2.2, and 3.2).

\textbf{Proposition 4.} \it Let a function $\varphi\in\mathcal{M}$ and positive real numbers $\varepsilon,\delta$ be given. Define a function $\psi\in\mathcal{B}$ by the formula
\begin{equation}\label{3f29}
\psi(t)=\begin{cases}
t^{\varepsilon/(\varepsilon+\delta)}\varphi(t^{1/(\varepsilon+\delta)})&  \text{if}\;\;t\geq 1\\
\varphi(1)&\text{if}\;\;0<t<1.
\end{cases}
\end{equation}
Then $\psi$ is an interpolation parameter, and
$$
\bigl[H^{s-\varepsilon}(G),H^{s+\delta}(G)\bigr]_{\psi}=H^{s,\varphi}(G)
\quad\mbox{for each}\quad s\in\mathbb{R}
$$
with equivalence of norms. Here, $G\in\{\mathbb{R}^{n},\Omega,\Gamma\}$.
If $G=\mathbb{R}^{n}$, then the equality of norms holds. \rm

We will also use the following three general properties of the interpolation.

\textbf{Proposition 5.} \it Let $X=[X_{0},X_{1}]$ and $Y=[Y_{0},Y_{1}]$ be admissible couples of Hilbert spaces, and let a linear mapping $T$ be given on $X_{0}$. Suppose that we have the bounded and Fredholm operators $T:X_{j}\rightarrow Y_{j}$, with $j=0,\,1$, that possess the common kernel and the common index. Then, for an arbitrary interpolation parameter $\psi\in\mathcal{B}$, the bounded operator $T:X_{\psi}\rightarrow Y_{\psi}$ is Fredholm, has the same kernel and the same index, and, moreover, its range $T(X_{\psi})=Y_{\psi}\cap T(X_{0})$. \rm

The proof of this proposition is given in \cite{MikhailetsMurach14} (Sect.~1.1.7).

\textbf{Proposition 6.} \it
Let $\bigl[X_{0}^{(j)},X_{1}^{(j)}\bigr]$, with $j=1,\ldots,r$, be a finite collection of admissible couples of Hilbert spaces. Then, for every function $\psi\in\mathcal{B}$, we have
$$
\biggl[\,\bigoplus_{j=1}^{r}X_{0}^{(j)},
\,\bigoplus_{j=1}^{r}X_{1}^{(j)}\biggr]_{\psi}=\,
\bigoplus_{j=1}^{r}\bigl[X_{0}^{(j)},X_{1}^{(j)}\bigr]_{\psi}
$$
with equality of norms. \rm

The proof of this proposition is given in \cite{MikhailetsMurach14} (Sect.~1.1.5).

The third property deals with the interpolation of subspaces connected with a linear bounded operator. Beforehand, we admit the following notation.

Let $H$, $\Phi$ and $\Psi$ be Hilbert spaces, with $\Phi\hookrightarrow\Psi$ continuously, and let a bounded linear  operator $T:H\rightarrow\Psi$ be given. We put
$$
(H)_{T,\Phi}:=\{u\in H:\,Tu\in\Phi\}.
$$
The linear space $(H)_{T,\Phi}$ is considered as a Hilbert space with respect to the graph norm
$$
\|u\|_{(H)_{T,\Phi}}:=\bigl(\|u\|_{H}^{2}+\|Tu\|_{\Phi}^{2}\bigr)^{1/2}.
$$
Note that this space does not depend on $\Psi$.

\textbf{Proposition 7.} \it Suppose that six separable Hilbert spaces $X_{0}$, $Y_{0}$, $Z_{0}$, $X_{1}$, $Y_{1}$, and $Z_{1}$ and three linear mappings $T$, $R$, and $S$ satisfy the following seven conditions:

\rm(i) \it the couples $X=[X_{0},X_{1}]$ and $Y=[Y_{0},Y_{1}]$ are admissible;

\rm(ii) \it the spaces $Z_{0}$ and $Z_{1}$ are subspaces of a certain linear space $E$;

\rm(iii) \it the continuous embedding $Y_{0}\hookrightarrow Z_{0}$ and $Y_{1}\hookrightarrow Z_{1}$ hold;

\rm(iv) \it the mapping $T$ is given on $X_{0}$ and defines the bounded operators $T:X_{0}\rightarrow Z_{0}$ and $T:X_{1}\rightarrow Z_{1}$;

\rm(v) \it the mapping $R$ is given on $E$ and defines the bounded operators $R:Z_{0}\rightarrow X_{0}$ and $R:Z_{1}\rightarrow X_{1}$;

\rm(vi) \it the mapping $S$ is given on $E$ and defines the bounded operators
$S:Z_{0}\rightarrow Y_{0}$ and $S:Z_{1}\rightarrow Y_{1}$;

\rm (vii) \it the equality $TR\,\omega=\omega+S\omega$ holds for each $\omega\in E$.

Then the couple of spaces $[\,(X_{0})_{T,Y_{0}},\,(X_{1})_{T,Y_{1}}\,]$ is admissible, and
\begin{equation*}
[\,(X_{0})_{T,Y_{0}},\,(X_{1})_{T,Y_{1}}\,]_{\psi} =(X_{\psi})_{T,Y_{\psi}}
\end{equation*}
up to equivalence of norms for an arbitrary interpolation parameter $\psi\in\mathcal{B}$. \rm

An analog of Proposition~7 appeared first in J.-L.~Lions and E.~Magenes' monograph \cite{LionsMagenes72} (Theorem 14.3), where the case of holomorphic interpolation was considered. This proposition is proved by V.~A.~Mikhailets and A.~A.~Murach in paper \cite{MikhailetsMurach06UMJ11} (Sect.~4); see also monograph \cite{MikhailetsMurach14} (Theorem 3.12).

\section{Proofs of the main results} 

Let us prove Theorems 1--3.

\textbf{\textit{Proof of Theorem}~1\textit{.}} We will first prove this theorem in the Sobolev case where $s\in\mathbb{Z}$, $s<2q$, and $\varphi \equiv1$. Consider the bounded and Fredholm operator \eqref{3f26} from Proposition~2. Certainly, the restriction of this operator to the space $H^{s,(2q)}_{A}(\Omega)$ is a bounded operator
\begin{equation}\label{3f30}
\Lambda:\,H_{A}^{s,(2q)}(\Omega)\oplus
\bigoplus_{k=1}^{\varkappa}H^{s+r_{k}-1/2}(\Gamma)\rightarrow
L_{2}(\Omega)\oplus\bigoplus_{j=1}^{q+\varkappa}H^{s-m_{j}-1/2}(\Gamma).
\end{equation}
Since $C^{\infty}(\overline{\Omega})$ is dense in $H_{A}^{s,(2q)}(\Omega)$, this operator is an extension by continuity of the mapping~\eqref{3f3}. Let $\Lambda_{s,(2q)}$ denote the operator \eqref{3f26}, and let $\Lambda_{s}$ denote the operator \eqref{3f30}. Evidently,  $\ker\Lambda_{s}=\ker\Lambda_{s,(2q)}$ and
\begin{equation*}
\mathrm{Ran}\,\Lambda_{s}=\bigl\{(f,g_{1},\ldots,g_{q+\varkappa})
\in\mathrm{Ran}\,\Lambda_{s,(2q)}:\,f\in L_{2}(\Omega)\bigr\},
\end{equation*}
with $\mathrm{Ran}$ denoting the range of the corresponding operator.
It follows from this by Proposition~2 that $\ker\Lambda_{s}=N$ and
\begin{equation*}
\mathrm{Ran}\,\Lambda_{s}=\biggl\{(f,g_{1},\ldots,g_{q+\varkappa})\in
L_{2}(\Omega)\oplus\bigoplus_{j=1}^{q+\varkappa}H^{s-m_{j}-1/2}(\Gamma):
\,\mbox{\eqref{3f12} is true}\biggr\}.
\end{equation*}
Therefore,
$$
\dim\ker\Lambda_{s}=\dim N<\infty\quad\mbox{and}\quad
\dim\mathrm{coker}\,\Lambda_{s}=\dim N^{+}<\infty.
$$
Thus, the bounded operator $\Lambda_{s}$ is Fredholm. In view of Proposition~3, we conclude that $\Lambda_{s}$ is the required operator \eqref{3f14} in Theorem~1, the set  $C^{\infty}(\overline{\Omega})$ being dense in $H^{s}_{A}(\Omega)$.

We will now prove Theorem~1 in the general situation with the help of interpolation with a function parameter between Sobolev spaces. Let $s\in\mathbb{R}$, $s<2q$, and $\varphi\in\mathcal{M}$. Choose an integer $p\geq2$ such that $s>-2q(p-1)$. The mapping \eqref{3f3} extends by continuity to the bounded and Fredholm operators
\begin{gather}\label{3f31}
\Lambda:\,H_{A}^{-2q(p-1)}(\Omega)\oplus
\bigoplus_{k=1}^{\varkappa}H^{-2q(p-1)+r_{k}-1/2}(\Gamma)\rightarrow
L_{2}(\Omega)\oplus
\bigoplus_{j=1}^{q+\varkappa}H^{-2q(p-1)-m_{j}-1/2}(\Gamma),\\
\Lambda:\,H^{2q}(\Omega)\oplus
\bigoplus_{k=1}^{\varkappa}H^{2q+r_{k}-1/2}(\Gamma)\rightarrow
L_{2}(\Omega)\oplus
\bigoplus_{j=1}^{q+\varkappa}H^{2q-m_{j}-1/2}(\Gamma).\label{3f32}
\end{gather}
This has just been proved for the first operator and is said in Proposition~2 for the second. These operators have the common kernel $N$ and the same index equaled to $\dim N-\dim N^{+}$. The range of the first operator satisfies the equation
\begin{equation}\label{3f33}
\begin{gathered}
\Lambda\bigl(H_{A}^{-2q(p-1)}(\Omega)\bigr)=\\=
\biggl\{(f,g_{1},\ldots,g_{q+\varkappa})\in
L_{2}(\Omega)\oplus
\bigoplus_{j=1}^{q+\varkappa}H^{-2q(p-1)-m_{j}-1/2}(\Gamma):
\mbox{\eqref{3f12} is true}\biggr\},
\end{gathered}
\end{equation}
and analogous formula holds for the second operator. Note that the second operator is a restriction of the first.

Let us define the interpolation parameter $\psi$ by formula \eqref{3f29} in which $\varepsilon:=s+2q(p-1)>0$ and $\delta:=2q-s>0$. Applying the interpolation with the function parameter $\psi$ to the couples of spaces in \eqref{3f31} and \eqref{3f32}, we obtain the bounded operator
\begin{equation}\label{3f34}
\begin{gathered}
\Lambda:\,\biggl[H_{A}^{-2q(p-1)}(\Omega)\oplus
\bigoplus_{k=1}^{\varkappa}H^{-2q(p-1)+r_{k}-1/2}(\Gamma),
H^{2q}(\Omega)\oplus
\bigoplus_{k=1}^{\varkappa}H^{2q+r_{k}-1/2}(\Gamma)\biggr]_{\psi}
\rightarrow\\
\rightarrow\biggl[L_{2}(\Omega)\oplus
\bigoplus_{j=1}^{q+\varkappa}H^{-2q(p-1)-m_{j}-1/2}(\Gamma),
L_{2}(\Omega)\oplus
\bigoplus_{j=1}^{q+\varkappa}H^{2q-m_{j}-1/2}(\Gamma)\biggr]_{\psi}.
\end{gathered}
\end{equation}
According to Proposition~5, this operator is Fredholm with kernel $N$ and index $\dim N-\dim N^{+}$.

Let us describe the interpolation spaces in \eqref{3f34}. Using Propositions 6 and 4 successively, we get
\begin{equation}\label{3f35}
\begin{gathered}
\biggl[H_{A}^{-2q(p-1)}(\Omega)\oplus
\bigoplus_{k=1}^{\varkappa}H^{-2q(p-1)+r_{k}-1/2}(\Gamma),
H^{2q}(\Omega)\oplus
\bigoplus_{k=1}^{\varkappa}H^{2q+r_{k}-1/2}(\Gamma)\biggr]_{\psi}=\\=
\bigl[H_{A}^{-2q(p-1)}(\Omega),H^{2q}(\Omega)\bigr]_{\psi}\oplus
\bigoplus_{k=1}^{\varkappa}
\bigl[H^{-2q(p-1)+r_{k}-1/2}(\Gamma),H^{2q+r_{k}-1/2}(\Gamma)\bigr]_{\psi}
=\\=\bigl[H_{A}^{-2q(p-1)}(\Omega),H^{2q}(\Omega)\bigr]_{\psi}\oplus
\bigoplus_{k=1}^{\varkappa}H^{s+r_{k}-1/2,\varphi}(\Gamma).
\end{gathered}
\end{equation}
Here, we take into account the equalities
\begin{equation}\label{3f36}
s-\varepsilon=-2q(p-1)\quad\mbox{and}\quad s+\delta=2q.
\end{equation}
Analogously,
\begin{equation}\label{3f37}
\begin{gathered}
\biggl[L_{2}(\Omega)\oplus
\bigoplus_{j=1}^{q+\varkappa}H^{-2q(p-1)-m_{j}-1/2}(\Gamma),
L_{2}(\Omega)\oplus
\bigoplus_{j=1}^{q+\varkappa}H^{2q-m_{j}-1/2}(\Gamma)\biggr]_{\psi}=\\=
L_{2}(\Omega)\oplus
\bigoplus_{j=1}^{q+\varkappa}H^{s-m_{j}-1/2,\varphi}(\Gamma).
\end{gathered}
\end{equation}
These equalities of spaces hold true up to equivalence of norms.

Let us prove that
\begin{equation}\label{3f38}
\bigl[H_{A}^{-2q(p-1)}(\Omega),H^{2q}(\Omega)\bigr]_{\psi}=
H_{A}^{s,\varphi}(\Omega)
\end{equation}
up to equivalence of norms. To this end, we use Proposition~7 in which
$$
X_{0}:=H^{-2q(p-1)}(\Omega),\quad X_{1}:=H^{2q}(\Omega),\quad
Y_{0}:=Y_{1}:=Z_{1}:=L_{2}(\Omega),\quad Z_{0}:=E:=H^{-2qp}(\Omega),
$$
and $T:=A$. Then
\begin{equation}\label{3f39}
H_{A}^{-2q(p-1)}(\Omega)=(X_{0})_{T,Y_{0}}\quad\mbox{and}\quad
H^{2q}(\Omega)=(X_{1})_{T,Y_{1}}.
\end{equation}
Note that the latter equality holds true up to equivalence of norms because $A$ is a bounded operator from $H^{\sigma}(\Omega)$ to $H^{\sigma-2q}(\Omega)$ for each $\sigma\in\mathbb{R}$, specifically, from $H^{2q}(\Omega)$ to $L_{2}(\Omega)$.

Evidently, conditions (i)--(iv) of Proposition~7 are fulfilled. We also need to give certain operators $R$ and $S$ satisfying conditions (v)--(vii). With this in mind, we use the known fact that the mapping $u\mapsto A^{p}A^{p+}u+u$ defines the isomorphism
\begin{equation}\label{3f40}
A^{p}A^{p+}+I:\,H^{\sigma}_{D}(\Omega)\leftrightarrow
H^{\sigma-4qp}(\Omega)\quad\mbox{for each}\quad \sigma\geq2qp
\end{equation}
(see, e.g., \cite{MikhailetsMurach14} (Lemma~3.1)). Here, of course, $A^{p}$ is the $p$-th iteration of $A$, then $A^{p+}$ is the formally adjoint operator to the differential operator $A^{p}$, and $I$ is the identity operator. Besides,
$$
H^{\sigma}_{D}(\Omega):=\bigl\{u\in H^{\sigma}(\Omega):\,
R_{\Gamma}D_{\nu}^{j-1}u=0\;\,\mbox{for each}\;\,j\in\{1,\ldots,2qp\}\bigr\}
$$
is a subspace of $H^{\sigma}(\Omega)$. The inverse of \eqref{3f40} defines the bounded linear operator
\begin{equation}\label{3f41}
(A^{p}A^{p+}+I)^{-1}:\,H^{\sigma}(\Omega)\rightarrow
H^{\sigma+4qp}(\Omega)\quad\mbox{for each}\quad \sigma\geq-2qp.
\end{equation}
We now let
$$
R:=A^{p-1}A^{p+}(A^{p}A^{p+}+I)^{-1}\quad\mbox{and}\quad S=-(A^{p}A^{p+}+I)^{-1}.
$$
Using \eqref{3f41}, we get the bounded operators
\begin{gather*}
R:\,Z_{0}=H^{-2qp}(\Omega)\rightarrow
H^{2qp-2qp-2q(p-1)}(\Omega)=X_{0},\\
R:\,Z_{1}=L_{2}(\Omega)\rightarrow H^{4qp-2qp-2q(p-1)}(\Omega)=X_{1},\\
S:\,Z_{0}=H^{-2qp}(\Omega)\rightarrow H^{2qp}(\Omega)\hookrightarrow
L_{2}(\Omega)=Y_{0},\\
S:\,Z_{1}=L_{2}(\Omega)\rightarrow H^{4qp}(\Omega)\hookrightarrow
L_{2}(\Omega)=Y_{1};
\end{gather*}
here, the embeddings are continuous. Besides,
$$
AR=AA^{p-1}A^{p+}(A^{p}A^{p+}+I)^{-1}=
(A^{p}A^{p+}+I-I)(A^{p}A^{p+}+I)^{-1}=I+S
$$
on $E=H^{-2qp}(\Omega)$. Thus, the rest conditions (v)--(vii) of Proposition~7 are satisfied.

Now, according to this proposition and in view of \eqref{3f39}, we get
\begin{equation*}
\bigl[H_{A}^{-2q(p-1)}(\Omega),H^{2q}(\Omega)\bigr]_{\psi}=
\bigl[(X_{0})_{T,Y_{0}},(X_{1})_{T,Y_{1}}\bigr]_{\psi}=
(X_{\psi})_{T,Y_{\psi}}.
\end{equation*}
Here, by Proposition~4 and \eqref{3f36}, we have
\begin{equation*}
X_{\psi}=\bigl[H^{-2q(p-1)}(\Omega),H^{2q}(\Omega)\bigr]_{\psi}=
H^{s,\varphi}(\Omega).
\end{equation*}
Besides,
$$
Y_{\psi}=\bigl[L_{2}(\Omega),L_{2}(\Omega)\bigr]_{\psi}=L_{2}(\Omega).
$$
These three formulas immediately gives the required equality \eqref{3f38}.

Note that
\begin{equation}\label{3f42}
\mbox{the set $C^{\infty}(\overline{\Omega})$ is dense in the space   $H^{s,\varphi}_{A}(\Omega)$}.
\end{equation}
It follows from the density of $C^{\infty}(\overline{\Omega})$ in $H^{2q}(\Omega)$ and from the dense continuous embedding of $H^{2q}(\Omega)$ in $H_{A}^{s,\varphi}(\Omega)$. The latter embedding is a general property of the interpolation used in~\eqref{3f38}.

It follows from the interpolation formulas \eqref{3f35}, \eqref{3f37}, and \eqref{3f38} that the bounded and Fredholm operator \eqref{3f34} acts between the spaces
\begin{equation}\label{3f43}
\Lambda:\,H^{s,\varphi}_{A}(\Omega)\oplus
\bigoplus_{k=1}^{\varkappa}H^{s+r_{k}-1/2,\varphi}(\Gamma)
\rightarrow L_{2}(\Omega)\oplus
\bigoplus_{j=1}^{q+\varkappa}H^{s-m_{j}-1/2,\varphi}(\Gamma).
\end{equation}
(Recall that this spaces are denoted by $\mathcal{D}^{s,\varphi}_{A}(\Omega,\Gamma)$ and $\mathcal{E}_{0,s,\varphi}(\Omega,\Gamma)$ respectively.)
We have proved that the operator \eqref{3f43} has the kernel $N$ and the index $\dim N-\dim N^{+}$. Furthermore, according to Proposition~5 and formula \eqref{3f33}, the range of this operator is equal to
$$
\mathcal{E}_{0,s,\varphi}(\Omega,\Gamma)\cap \Lambda\bigl(H_{A}^{-2q(p-1)}(\Omega)\bigr)=
\bigl\{(f,g_{1},\ldots,g_{q+\varkappa})\in
\mathcal{E}_{0,s,\varphi}(\Omega,\Gamma):
\mbox{\eqref{3f12} is true}\bigr\}.
$$
Finally, in view of \eqref{3f42}, the operator \eqref{3f43} is an extension by continuity of the mapping \eqref{3f3}. Thus, this operator is the required operator \eqref{3f14}.

Theorem~1 is proved.

\textbf{\textit{Proof of Theorem}~2\textit{.}} It follows from the isomorphism \eqref{3f-isom} that
\begin{equation}\label{3f-2a}
\inf\bigl\{\,\|(u,v)+(u^{(0)},v^{(0)})\|_
{\mathcal{D}^{s,\varphi}_{A}(\Omega,\Gamma)}:
\,(u^{(0)},v^{(0)})\in N\,\bigr\}\leq
c_{0}\,\|\Lambda(u,v)\|_{\mathcal{E}_{0,s,\varphi}(\Omega,\Gamma)}
\end{equation}
for each $(u,v)\in\mathcal{D}^{s,\varphi}_{A}(\Omega,\Gamma)$, with  $c_{0}$ being the norm of the inverse operator to \eqref{3f-isom}. Since $N$ is a finite-dimensional subspace of both spaces $\mathcal{D}^{s,\varphi}_{A}(\Omega,\Gamma)$ and $\mathcal{D}_A^{s-\sigma,\varphi}(\Omega,\Gamma)$, the norms in them are equivalent on~$N$. Hence, for each $(u^{(0)},v^{(0)})\in N$, we can write \begin{equation*}
\|(u^{(0)},v^{(0)})\|_{\mathcal{D}^{s,\varphi}_{A}(\Omega,\Gamma)}
\leq c_{1}
\|(u^{(0)},v^{(0)})\|_{\mathcal{D}^{s-\sigma,\varphi}_{A}(\Omega,\Gamma)}
\end{equation*}
with some number $c_{1}>0$, which is independent of $(u,v)$ and $(u^{(0)},v^{(0)})$. Besides,
\begin{gather*}
\|(u^{(0)},v^{(0)})\|_{\mathcal{D}^{s-\sigma,\varphi}_{A}(\Omega,\Gamma)}
\leq\|(u,v)+(u^{(0)},v^{(0)})\|_
{\mathcal{D}^{s-\sigma,\varphi}_{A}(\Omega,\Gamma)}+
\|(u,v)\|_{\mathcal{D}^{s-\sigma,\varphi}_{A}(\Omega,\Gamma)}
\leq\\ \leq c_{2}\,
\|(u,v)+(u^{(0)},v^{(0)})\|_{\mathcal{D}^{s,\varphi}_{A}(\Omega,\Gamma)}+
\|(u,v)\|_{\mathcal{D}^{s-\sigma,\varphi}_{A}(\Omega,\Gamma)}.
\end{gather*}
Here, $c_{2}$ is the norm of the continuous embedding operator $\mathcal{D}^{s,\varphi}_{A}(\Omega,\Gamma)\hookrightarrow
\mathcal{D}^{s-\sigma,\varphi}_{A}(\Omega,\Gamma)$. Hence,
\begin{gather*}
\|(u,v)\|_{\mathcal{D}^{s,\varphi}_{A}(\Omega,\Gamma)}\leq
\|(u,v)+(u^{(0)},v^{(0)})\|_{\mathcal{D}^{s,\varphi}_{A}(\Omega,\Gamma)}+
\|(u^{(0)},v^{(0)})\|_{\mathcal{D}^{s,\varphi}_{A}(\Omega,\Gamma)}\leq\\
\leq
\|(u,v)+(u^{(0)},v^{(0)})\|_{\mathcal{D}^{s,\varphi}_{A}(\Omega,\Gamma)}+
c_{1}
\|(u^{(0)},v^{(0)})\|_{\mathcal{D}^{s-\sigma,\varphi}_{A}(\Omega,\Gamma)}
\leq\\ \leq(1+c_{1}c_{2})\,
\|(u,v)+(u^{(0)},v^{(0)})\|_{\mathcal{D}^{s,\varphi}_{A}(\Omega,\Gamma)}+
c_{1}\,\|(u,v)\|_{\mathcal{D}^{s-\sigma,\varphi}_{A}(\Omega,\Gamma)}.
\end{gather*}
Passing in this inequality to the infimum over all $(u^{(0)},v^{(0)})\in N$ and using the bound \eqref{3f-2a}, we arrive at the estimate \eqref{3f20}, namely, we get
$$
\|(u,v)\|_{\mathcal{D}^{s,\varphi}_{A}(\Omega,\Gamma)}\leq
(1+c_{1}c_{2})\,c_{0}\,
\|\Lambda(u,v)\|_{\mathcal{E}_{0,s,\varphi}(\Omega,\Gamma)}
+ c_{1}\,\|(u,v)\|_{\mathcal{D}_A^{s-\sigma,\varphi}(\Omega,\Gamma)}\bigl).
$$
Theorem~2 is proved.

\textbf{\textit{Proof of Theorem}~3\textit{.}} We first consider the case of $s<2q$. By virtue of Theorem~1, the vector $(f,g)=\Lambda(u,v)$ satisfies \eqref{3f12}. Note that $(f,g)\in\mathcal{E}_{0,s,\varphi}(\Omega,\Gamma)$ by the condition \eqref{3f22}. Therefore, according to Theorem~1, we get the inclusion $(f,g)\in\Lambda(\mathcal{D}^{s,\varphi}_{A}(\Omega,\Gamma))$.
Thus, together with $\Lambda(u,v)=(f,g)$, we have the equality   $\Lambda(u',v')=(f,g)$ for a certain vector    $(u',v')\in\mathcal{D}^{s,\varphi}_{A}(\Omega,\Gamma)$. Hence,  $\Lambda(u-u',v-v')=0$, which implies that
$$
(u-u',v-v')\in N\subset C^{\infty}(\overline{\Omega})\times
\bigl(C^{\infty}(\Gamma)\bigr)^{\varkappa}.
$$
Thus,
$$
(u,v)=(u',v')+(u-u',v-v')\in\mathcal{D}^{s,\varphi}_{A}(\Omega,\Gamma).
$$
Theorem~3 has been proved in the case of $s<2q$.

Let us examine the case of $s\geq2q$. According to what has just been proved, the condition
$$
g_{j}\in H^{s-m_{j}-1/2,\varphi}(\Gamma)\subset H^{2q-1/4-m_{j}-1/2,\varphi}(\Gamma)
\quad\mbox{for each}\quad j\in\{1,...,q+\varkappa\}
$$
implies the inclusion
\begin{equation*}
(u,v)\in H^{2q-1/4,\varphi}(\Omega)\oplus
\bigoplus_{k=1}^{\varkappa}H^{2q-1/4+r_{k}-1/2,\varphi}(\Gamma).
\end{equation*}
Repeating the above reasoning and using Proposition~1 instead of Theorem~1, we deduce from this inclusion and the condition of Theorem~3, that the vector $(u,v)$ satisfies \eqref{3f-a} and \eqref{3f-b}.

Theorem~3 is proved.

\end{document}